\documentstyle[12pt,amssymb,amsfonts,twoside]{article}

\begin{titlepage}

\title{Transgression of the index gerbe}
\author{Ulrich Bunke\thanks{G\"ottingen, bunke@uni-math.gwdg.de} }
\end{titlepage}
\begin{document}
\maketitle
\begin{abstract}
The holonomy of an unitary line bundle with connection over some base space
$B$ is a $U(1)$-valued function on the loop space $LB$. In a parallel
manner, the holonomy of a gerbe with connection on $B$ is a line bundle with
connection over $LB$.

Given a family of graded Dirac operators on $B$ and some additional
geometric data one can define the determinant line bundle with Quillen metric
and Bismut-Freed connection. According to Witten, Bismut-Freed the holonomy of 
this determinant bundle can be expressed in terms of an adiabatic limit of eta
invariants of an associated family of Dirac operators over $LB$.

Recently, for a family of ungraded Dirac operators on $B$  J. Lott constructed
an index gerbe with connection. In the present paper we show, in analogy to the
holonomy formula for the determinant bundle, that the holonomy of the index
gerbe coincides with an adiabatic limit of determinant bundles of the
associated  family of Dirac operators over $LB$.

\end{abstract}

\maketitle
% Definitionen
\newcommand{\LIM}{{\rm LIM}}
\newcommand{\diag}{{\rm diag}}
\newcommand{\proof}{{\it Proof.$\:\:\:\:$}}
 \newcommand{\dist}{{\rm dist}}
\newcommand{\kaaa}{{\frak k}}
\newcommand{\paaa}{{\frak p}}
\newcommand{\vp}{{\varphi}}
\newcommand{\taaa}{{\frak t}}
\newcommand{\haaa}{{\frak h}}
\newcommand{\R}{{\Bbb R}}
\newcommand{\Hh}{{\bf H}}
\newcommand{\Rep}{{\rm Rep}}
\newcommand{\Hb}{{\Bbb H}}
\newcommand{\Q}{{\Bbb Q}}
\newcommand{\str}{{\rm str}}
\newcommand{\Ind}{{\rm ind}}
\newcommand{\triv}{{\rm triv}}
\newcommand{\Z}{{\Bbb Z}}
\newcommand{\bD}{{\bf D}}
\newcommand{\bF}{{\bf F}}
\newcommand{\tX}{{\tt X}}
\newcommand{\Cliff}{{\rm Cliff}}
\newcommand{\tY}{{\tt Y}}
\newcommand{\tZ}{{\tt Z}}
\newcommand{\tV}{{\tt V}}
\newcommand{\tR}{{\tt R}}
\newcommand{\Fam}{{\rm Fam}}
\newcommand{\Cusp}{{\rm Cusp}}
\newcommand{\bT}{{\bf T}}
\newcommand{\bK}{{\bf K}}
\newcommand{\K}{{\Bbb K}}
\newcommand{\tH}{{\tt H}}
\newcommand{\bS}{{\bf S}}
\newcommand{\bB}{{\bf B}}
\newcommand{\tW}{{\tt W}}
\newcommand{\tF}{{\tt F}}
\newcommand{\bA}{{\bf A}}
\newcommand{\bL}{{\bf L}}
 \newcommand{\bom}{{\bf \Omega}}
\newcommand{\bundle}{{bundle}}
\newcommand{\ch}{{\bf ch}}
\newcommand{\ve}{{\varepsilon}}
\newcommand{\C}{{\Bbb C}}
\newcommand{\gen}{{\rm gen}}
\newcommand{\cTop}{{{\cal T}op}}
\newcommand{\bP}{{\bf P}}
\newcommand{\Naaa}{{\bf N}}
\newcommand{\image}{{\rm image}}
\newcommand{\gaaa}{{\frak g}}
\newcommand{\zaaa}{{\frak z}}
\newcommand{\saaa}{{\frak s}}
\newcommand{\laaa}{{\frak l}}
\newcommand{\stimes}{{\times\hspace{-1mm}\bf |}}
\newcommand{\ausg}{{\rm end}}
\newcommand{\bff}{{\bf f}}
\newcommand{\maaa}{{\frak m}}
\newcommand{\aaaa}{{\frak a}}
\newcommand{\naaa}{{\frak n}}
\newcommand{\brr}{{\bf r}}
\newcommand{\res}{{\rm res}}
\newcommand{\Aut}{{\rm Aut}}
\newcommand{\Pol}{{\rm Pol}}
\newcommand{\Tr}{{\rm Tr}}
\newcommand{\cT}{{\cal T}}
\newcommand{\dom}{{\rm dom}}
\newcommand{\Line}{{\rm Line}}
\newcommand{\db}{{\bar{\partial}}}
\newcommand{\g}{{\gaaa}}
\newcommand{\cZ}{{\cal Z}}
\newcommand{\cH}{{\cal H}}
\newcommand{\cM}{{\cal M}}
\newcommand{\interi}{{\rm int}}
\newcommand{\singsupp}{{\rm singsupp}}
\newcommand{\cE}{{\cal E}}
\newcommand{\ccR}{{\cal R}}
\newcommand{\hol}{{\rm hol}}
\newcommand{\cV}{{\cal V}}
\newcommand{\cY}{{\cal Y}}
\newcommand{\cW}{{\cal W}}
\newcommand{\cI}{{\cal I}}
\newcommand{\cC}{{\cal C}}
\newcommand{\mod}{{\rm mod}}
\newcommand{\cK}{{\cal K}}
\newcommand{\cA}{{\cal A}}
\newcommand{\cEp}{{{\cal E}^\prime}}
\newcommand{\cU}{{\cal U}}
\newcommand{\Hom}{{\mbox{\rm Hom}}}
\newcommand{\vol}{{\rm vol}}
\newcommand{\cO}{{\cal O}}
\newcommand{\End}{{\mbox{\rm End}}}
\newcommand{\Ext}{{\mbox{\rm Ext}}}
\newcommand{\rk}{{\rm rank}}
\newcommand{\im}{{\mbox{\rm im}}}
\newcommand{\sign}{{\rm sign}}
\newcommand{\spann}{{\mbox{\rm span}}}
\newcommand{\symm}{{\mbox{\rm symm}}}
\newcommand{\cF}{{\cal F}}
\newcommand{\cD}{{\cal D}}
\newcommand{\Ree}{{\rm Re }}
\newcommand{\Res}{{\mbox{\rm Res}}}
\newcommand{\Imm}{{\rm Im}}
\newcommand{\inter}{{\rm int}}
\newcommand{\clo}{{\rm clo}}
\newcommand{\tg}{{\rm tg}}
\newcommand{\ee}{{\rm e}}
\newcommand{\Li}{{\rm Li}}
\newcommand{\cN}{{\cal N}}
 \newcommand{\conv}{{\rm conv}}
\newcommand{\op}{{\mbox{\rm Op}}}
\newcommand{\tr}{{\mbox{\rm tr}}}
\newcommand{\cs}{{c_\sigma}}
\newcommand{\ctg}{{\rm ctg}}
\newcommand{\degg}{{\mbox{\rm deg}}}
\newcommand{\Ad}{{\mbox{\rm Ad}}}
\newcommand{\ad}{{\mbox{\rm ad}}}
\newcommand{\codim}{{\rm codim}}
\newcommand{\Gr}{{\mathrm{Gr}}}
\newcommand{\coker}{{\rm coker}}
\newcommand{\id}{{\mbox{\rm id}}}
\newcommand{\ord}{{\rm ord}}
\newcommand{\nat}{{\Bbb  N}}
\newcommand{\supp}{{\rm supp}}
\newcommand{\sing}{{\mbox{\rm sing}}}
\newcommand{\spec}{{\mbox{\rm spec}}}
\newcommand{\Ann}{{\mbox{\rm Ann}}}
\newcommand{\aca}{{\aaaa_\C^\ast}}
\newcommand{\acag}{{\aaaa_{\C,good}^\ast}}
\newcommand{\acage}{{\aaaa_{\C,good}^{\ast,extended}}}
\newcommand{\tck}{{\tilde{\ck}}}
\newcommand{\tnk}{{\tilde{\ck}_0}}
\newcommand{\ceep}{{{\cal E}(E)^\prime}}
 \newcommand{\ncE}{{{}^\naaa\cE}}
 \newcommand{\Or}{{\rm Or }}
\newcommand{\Diff}{{\cal D}iff}
\newcommand{\cB}{{\cal B}}
\newcommand{\hc}{{{\cal HC}(\gaaa,K)}}
\newcommand{\hcma}{{{\cal HC}(\maaa_P\oplus\aaaa_P,K_P)}}
\def\imath{{\rm i}}
\newcommand{\vsl}{{V_{\sigma_\lambda}}}
\newcommand{\czg}{{\cZ(\gaaa)}}
\newcommand{\csl}{{\chi_{\sigma,\lambda}}}
\newcommand{\cR}{{R}}
\def\hB{\hspace*{\fill}$\Box$ \newline\noindent}
\newcommand{\varho}{\varrho}
\newcommand{\ind}{{\rm index}}
\newcommand{\Indu}{{\rm Ind}}
\newcommand{\Fin}{{\mbox{\rm Fin}}}
\newcommand{\cS}{{S}}
\newcommand{\orig}{{\cal O}}
\def\hB{\hspace*{\fill}$\Box$ \\[0.5cm]\noindent}
\newcommand{\cL}{{\cal L}}
 \newcommand{\cG}{{\cal G}}
\newcommand{\npci}{{\naaa_P\hspace{-1.5mm}-\hspace{-2mm}\mbox{\rm coinv}}}
\newcommand{\pki}{{(\paaa,K_P)\hspace{-1.5mm}-\hspace{-2mm}\mbox{\rm inv}}}
\newcommand{\mki}{{(\maaa_P\oplus \aaaa_P, K_P)\hspace{-1.5mm}-\hspace{-2mm}\mbox{\rm inv}}}
\newcommand{\Mat}{{\rm Mat}}
\newcommand{\npi}{{\naaa_P\hspace{-1.5mm}-\hspace{-2mm}\mbox{\rm inv}}}
\newcommand{\ngp}{{N_\Gamma(\pi)}}
\newcommand{\gbg}{{\Gamma\backslash G}}
\newcommand{\gkm}{{ Mod(\gaaa,K) }}
\newcommand{\ggkm}{{  (\gaaa,K) }}
\newcommand{\pkm}{{ Mod(\paaa,K_P)}}
\newcommand{\ppkm}{{  (\paaa,K_P)}}
\newcommand{\makm}{{Mod(\maaa_P\oplus\aaaa_P,K_P)}}
\newcommand{\mmakm}{{ (\maaa_P\oplus\aaaa_P,K_P)}}
\newcommand{\cP}{{\cal P}}
\newcommand{\gm}{{Mod(G)}}
\newcommand{\gk}{{\Gamma_K}}
\newcommand{\La}{{\cal L}}
\newcommand{\cug}{{\cU(\gaaa)}}
\newcommand{\cuk}{{\cU(\kaaa)}}
\newcommand{\dc}{{C^{-\infty}_c(G) }}
\newcommand{\gdk}{{\gaaa/\kaaa}}
\newcommand{\dgkm}{{ D^+(\gaaa,K)-\mbox{\rm mod}}}
\newcommand{\dgm}{{D^+G-\mbox{\rm mod}}}
\newcommand{\vect}{{\C-\mbox{\rm vect}}}
 \newcommand{\cig}{{C^{ \infty}(G)_{K} }}
\newcommand{\gami}{{\Gamma\hspace{-1.5mm}-\hspace{-2mm}\mbox{\rm inv}}}
\newcommand{\cQ}{{\cal Q}}
\newcommand{\mmap}{{Mod(M_PA_P)}}
\newcommand{\bbbz}{{\bf Z}}
 \newcommand{\cX}{{\cal X}}
\newcommand{\bH}{{\bf H}}
\newcommand{\pr}{{\rm pr}}
\newcommand{\bX}{{\bf X}}
\newcommand{\bY}{{\bf Y}}
\newcommand{\bZ}{{\bf Z}}
\newcommand{\bV}{{\bf V}}
\newcommand{\Gerbe}{{\rm Gerbe}}
\newcommand{\gerbe}{{\rm gerbe}}
\newcommand{\hA}{{\bf \hat A}}

\newtheorem{prop}{Proposition}[section]
\newtheorem{lem}[prop]{Lemma}
\newtheorem{ddd}[prop]{Definition}
\newtheorem{theorem}[prop]{Theorem}
\newtheorem{kor}[prop]{Corollary}
\newtheorem{ass}[prop]{Assumption}
\newtheorem{con}[prop]{Conjecture}
\newtheorem{prob}[prop]{Problem}
\newtheorem{fact}[prop]{Fact}

\tableofcontents

\parskip3ex
\section{Families and loops}

\subsection{Geometric families}\label{geo}

In this subsection we combine various geometric structures on a smooth fibre
bundle by introducing the notion of a geometric family.

Let $\pi:M\rightarrow B$ be a smooth fibre bundle with closed fibres.
The base $B$ of this bundle may be infinite dimensional, and the fibres of
$\pi$ may consist of several components of different dimensions. By $T^v\pi$ we
denote the vertical bundle.  We assume that $T^v\pi$ is oriented and carries a
spin structure. Furthermore, let $g^{T^v\pi}$ be a vertical Riemannian metric,
and let $T^h\pi$ be a horizontal distribution. Finally, let
$(V,h^V,\nabla^V)$ be a hermitean vector bundle with metric connection over
$M$.

\begin{ddd}
A geometric family over $B$ 
is a collection of objects $\pi:M\rightarrow B$, spin structure and
orientation of $T^v\pi$, $g^{T^v\pi}$, $T^h\pi$, and $(V,h^V,\nabla^V)$
as introduced above.
\end{ddd}

We now describe some natural operations with geometric families.
If $f:B^\prime\rightarrow B$ is a smooth map and $\cE$ is a geometric family
over $B$, then we define the geometric family $f^*\cE$ over $B^\prime$ as
the collection of the following objects.
As the smooth fibre bundle we take $f^*\pi:f^*M\rightarrow B^\prime$.
There is a natural map $F:f^*M\rightarrow M$ which is a diffeomorphism
fibrewise. In particular, we have an isomorphism
$dF_{|T^vf^*\pi}:T^vf^*\pi\stackrel{\sim}{\rightarrow} F^* T^v\pi$.
Using this isomorphism we obtain the  induced orientation and spin structure 
on $T^vf^*\pi$. We define the vertical metric $g^{T^vf^*\pi}$ such that
$dF_{|T^vf^*\pi}$ becomes an isometry.
The induced horizontal distribution $T^hf^*\pi$ is defined as the kernel
of the composition $\pr_{F^*T^v\pi}\circ dF:Tf^*M\rightarrow F^*T^v\pi$,
where $\pr_{F^*T^v\pi}:F^*TM\rightarrow F^*T^v\pi$ is the projection along
$F^*T^h\pi$. The hermitean bundle with connection over $f^*M$ is the pull back
$F^*V$ with induced metric $\nabla^{F^*V}$ and connection $h^{F^*V}$.

If $\cE_0,\cE_1$ are geometric families over $B_0, B_1$, then we can form the
union $\cE_0\cup \cE_1$ in a natural way. It is a geometric family over
$B_0\cup B_1$. 
In particular, if $B=B_0=B_1$, then we can compose the projection
$M_0\cup M_1\rightarrow B\cup B$ with the natural covering map $B\cup
B\rightarrow B$. The resulting geometric family over $B$ is the relative
union $\cE_0\cup_B\cE_1$.
The fibre of $\cE_0\cup_B \cE_1$ is the disjoint
union of the fibres of $\cE_0$ and $\cE_1$.

We can also form the product $\cE:=\cE_0\times \cE_1$, a geometric family over
$B_0\times B_1$. On the level of the hermitean bundles we take here
$V:=\pr_{M_0}^*V_0\otimes \pr^*_{M_1}V_1$, where $\pr_{M_i}:
M_0\times M_1\rightarrow M_i$ are the natural projections.
Again, if $B=B_0=B_1$, then we can form the pull-back of $\cE$ with respect to
the diagonal embedding $B\rightarrow B\times B$. The resulting geometric family
is the fibre product $\cE_0\times_B\cE_1$.

The main construction for the purpose of the present paper is the loop $L\cE$
of a geometric family $\cE$ over $B$. Here $L\cE$ is a geometric family over
$LB\times \R_+$, where $LB$ is the free loop space of $B$.

Consider the evaluation map $\tilde{ev}:LB\times \R_+\times S^1\rightarrow B$
given by $\tilde{ev}(\gamma,\epsilon,u)=\gamma(u)$.
First we consider the geometric family $\tilde{ev}^*\cE$ over $LB\times \R_+$.
The bundle of $L\cE$ will be the composition
$L\pi:\tilde{ev}^*M\rightarrow LB\times \R_+\times
S^1\stackrel{\pr}{\rightarrow} LB\times \R_+$.
Thus the fibre of $L\pi$ over $(\gamma,\epsilon)$ is the total space
of the pull-back bundle $\gamma^*M\rightarrow S^1$.

The differential $d\tilde{ev}^*\pi$ induces an isomorphism
$\Phi:T^h\tilde{ev}^*\pi\stackrel{\sim}{\rightarrow} (\tilde{ev}^*\pi)^*
(TLB\oplus T\R_+\oplus TS^1)$. Then we have $T^vL\pi\cong
T^v\tilde{ev}^*\pi\oplus \Phi^{-1} (\tilde{ev}^*\pi)^* TS^1$, and we define
$T^hL\pi:=\Phi^{-1}(TLB\oplus T\R_+)$.
We equip $TS^1$
with the bounding spin structure and the orientation obtained from the
identification $S^1:=\R/\Z$. Together with the orientation and spin structure
of  $T^v\tilde{ev}^*\pi$ this induces an orientation and a spin structure on 
$T^vL\pi$. The description of $S^1$ above also induces a metric $g^{TS^1}$ of
volume one. We define the vertical metric of $L\cE$ by 
$g^{T^v\tilde{ev}^*\pi}\oplus \frac{1}{\epsilon^2} \Phi^* (\tilde{ev}^*\pi)^*
g^{TS^1}$. The hermitean vector bundle with connection of $L\cE$ is
$(\tilde{Ev}^*V,h^{\tilde{Ev}^*V},\nabla^{\tilde{Ev}^*V})$,
where $\tilde{Ev}:\tilde{ev}^*M\rightarrow M$ is the natural map.

Let $I_\epsilon:LB\rightarrow LB\times \R_+$ be given by
$I_\epsilon(\gamma):=(\gamma,\epsilon)$.
Then we set $L_\epsilon \cE:=I_\epsilon^* L\cE$.

Note that $L(\cE_0\cup_B\cE_1)=L\cE_0\cup_{LB\times \R_+} L\cE_1$.

If $\pi:M\rightarrow B$ is a smooth fibre bundle and $g^M$ is a Riemannian
metric on $M$, then it induces a Riemannian metric $g^{T^v\pi}$ by
restriction, and it determines a horizontal distribution $T^h\pi$ as the
orthogonal complement of $T^v\pi$. 

On the other hand, if $\pi:M\rightarrow B$ is a smooth fibre bundle  equiped
with a vertical metric $g^{T^v\pi}$ and a horizontal distribution $T^h\pi$,
then any Riemannian metric $g^B$ on $B$ gives rise to a family of Riemannian
metrics $g_\epsilon^M:=g^{T^v\pi}\oplus \frac{1}{\epsilon^2} \pi^* g^B$,
$\epsilon\in \R_+$. The definition  of the fibre wise metric of the loop of a
geometric family is a special case of this construction.
The limit $\epsilon\to 0$ is called the adiabatic limit. In the present paper
we will meet adiabatic limits of various geometric (e.g. curvature)  and 
spectral geometric (e.g. $\eta$-invariant) quantities associated to
the Riemannian manifold $(M,g_\epsilon^M)$.

Given $\pi:M\rightarrow B$, $g^{T^v\pi}$ and  $T^h\pi$, we construct a
connection $\nabla^{T^v\pi}$ as follows. We choose any Riemannian metric
$g^B$ and let $\nabla^{T^v\pi}$ be the projection of the Levi-Civita
connection of $g^M_1$ to $T^v\pi$. The connection $\nabla^{T^v\pi}$ is
independent of the choice of $g^B$ (comp. \cite{berlinegetzlervergne92}, Proposition 
10.2)

\subsection{Families of Dirac operators}\label{dirac}

Given a geometric family $\cE$ over $B$ we have a family of twisted
Dirac operators $D(\cE):=(D_b(\cE))_{b\in B}$.
By $\Gamma(\cE)$ we denote the infinite dimensional bundle over $B$ the fibre over $b\in B$
of which  is the space of sections $C^\infty(M_b,S(M_b)\otimes
V_{|M_b})$, where $M_b:=\pi^{-1}(b)$ and $S(M_b)$ is the spinor bundle of
$M_b$. We can consider $D(\cE)$ as a section of $\End(\Gamma(\cE))$.
The bundle $\Gamma(\cE)$ carries a hermitean scalar product and a natural
hermitean connection $\nabla^{\Gamma(\cE)}$ (see
\cite{berlinegetzlervergne92}, Proposition  9.13). Let $T\in
C^\infty(M,\Hom(\Lambda^2T^h\pi,T^v\pi))$ be the curvature tensor 
(it is the negative of the tensor $\Omega$ defined in
\cite{berlinegetzlervergne92}, p.321) of the
horizontal distribution. We form $c(T)\in \Omega^2(B,\End(\Gamma(\cE)))$ such
that $c(T)(X,Y)$ is given by Clifford multiplication by $T(X^h,Y^h)$,
where $X^h,Y^h$ denote the horizontal lifts of $X,Y$.

So a geometric bundle gives rise to the scaled Bismut super connection (see
\cite{berlinegetzlervergne92}, Proposition  10.15) 
$$A_s(\cE):= s D +
\nabla^{\Gamma(\cE)} + \frac{1}{4s} c(T)$$ if the fibres of $\pi$ are even
dimensional. If they are odd dimensional, then we set 
$$A_s(\cE):= s \sigma D + \nabla^{\Gamma(\cE)} + \frac{\sigma}{4s} c(T) \ ,$$
where $\sigma$ is an additional odd variable satisfying $\sigma^2=1$.

If the base is finite-dimensional,  $TB$ is oriented and has a
spin structure, and $g^B$ is a Riemannian metric on the base, then
$TM=T^v\pi\oplus \pi^*TB$ has an induced orientation and spin structure.
By $D_\epsilon$ we denote the $V$-twisted Dirac operator on $M$ associated to
the Riemannian metric $g_\epsilon^M$ introduced in Subsection \ref{geo}.

\subsection{Transgression in $K$-theory and cohomology}

Given a geometric family $\cE$ over a base $B$ there is a close relation
between the index  $\ind(\cE)\in K^*(B)$ of the family $D(\cE)$ and
$\ind(L_1\cE)\in K^{*+1}(LB)$ of the family $D(L_1\cE)$.

The choice of the spin structure of $TS^1$ provides a $K$-orientation
of the fibres of $\pr : LB\times S^1\rightarrow LB$ and therefore a map
$\pr_!:K^*(LB\times S^1)\rightarrow K^{*-1}(LB)$.
We define the transgression map
$T: K^*(B)\rightarrow K^{*-1}(LB)$ by $T:=\pr_!\circ ev^*$.
Then we have
\begin{equation}
\ind(L_1\cE)= T(\ind(\cE))\ .
\end{equation}
In \cite{daizhang96} this result is attributed to Wu \cite{wu97}. 
The orientation of $S^1$ induces an integration-over-the-fibre
$\pr_!:H^*(LB\times S^1,\R)\rightarrow H^{*-1}(B,\R)$ and therefore a
transgression map in cohomology $T:= \pi_!\circ ev^*:H^*(LB\times
S^1,\R)\rightarrow H^{*-1}(LB,\R)$. Transgression is compatible with the Chern
character, i.e. $T\circ \ch =\ch\circ T$.
In particular,
$$T\circ \ch(\ind(\cE))=\ch(\ind(L_1\cE))\ .$$
The goal of the present paper is a refinement of this
equality replacing the Chern classes by their Delinge cohomology
valued refinements due to Cheeger-Simons \cite{cheegersimons83}.
To be honest, the corresponding equality for the component in first Deligne
cohomology  is known for a while, and it is equivalent to the holonomy
formula for the determinant line bundle (Theorem \ref{bf}). In the present
paper we consider the degree two component.

\subsection{Transgression of line bundles with connection}

Let $\Line(B)$ denote the set of isomorphism classes of hermitean line bundles
with connection over $B$. Note that the tensor product induces a group
structure on  $\Line(B)$. The inverse is given by the dual bundle.

We define a transgression $T:\Line(B)\rightarrow C^\infty(LB,U(1))$
as follows. Given $\cL:=(L,h^L,\nabla^L)$ and $\gamma\in L(B)$, then let
$\hol_\gamma(\cL)$ be the holonomy of $\cL$ along the loop $\gamma$.
We define  $T(\cL)(\gamma):=\hol_\gamma(\cL)$. 

The transgression is a group homomorphism. Note that $\cL$ is completely
determined by $T(\cL)$. Furthermore, $T(\cL)$ is invariant under the natural
action of $\Diff(S^1)$ on $LB$.

\subsection{Transgression of the determinant line bundle}

Let $\cE$ be a geometric family over $B$ with even dimensional fibres.
Then we have a determinant bundle of $D(\cE)$ with hermitean metric (Quillen
metric) and hermitean connection (Bismut-Freed connection) which we denote by
$\det(\cE)\in \Line(B)$ (see \cite{bismutfreed861}, or
\cite{berlinegetzlervergne92}, Ch. 9.7) . 

If $D$ is a twisted Dirac operator on a closed oriented odd-dimensional 
Riemannian spin manifold, then  
following Dai/Freed 
\cite{daifreed94}, we define 
$$\tau(D):=\exp(2\pi\imath \frac{\eta(D)+\dim\ker
D}{2})\ ,$$ 
where $\eta(D):=\frac{1}{\pi^{1/2}}\int_0^\infty \Tr
D\ee^{-uD^2}\frac{du}{u^{1/2}}$ is the $\eta$-invariant of $D$ first
introduced by Atiyah-Patodi-Singer \cite{atiyahpatodisinger75}.  Applying this
to the family $D(L_\epsilon\cE)$ we obtain a  function
$\tau(D(L_\epsilon\cE))\in C^\infty(LB,U(1))$. The adiabatic limit
$\lim_{\epsilon\to 0}\tau(D(L_\epsilon\cE))$ exists locally uniformly. The main result of
Bismut-Freed \cite{bismutfreed861}, \cite{bismutfreed862} (see also
\cite{daifreed94} for a different proof) is \begin{theorem}\label{bf}
$$T\det(\cE) = \lim_{\epsilon\to 0}\tau(D(L_\epsilon\cE))$$
\end{theorem}

\subsection{Transgression of gerbes with connection}

Let $\Gerbe(B)$ denote the set of isomorphism classes of gerbes with
connection over $B$ (see Subsection \ref{gerbe} for a definition).
There is a product of gerbes making $\Gerbe(B)$ into an abelian group.

In Subsection \ref{trans} we will define a transgression map 
$T:\Gerbe(B)\rightarrow \Line(LB)$.

There is a characteristic class $v:\Gerbe(B)\rightarrow H^3(B,\Z)$ classifying
isomorphism classes of gerbes (forgetting the connection).
\begin{lem}\label{com1}
 If
$\cG\in\Gerbe(B)$,
then $-T(v(\cG))=c_1(T(\cG))$.
\end{lem}
Here $c_1(T(\cG))\in H^2(LB,\Z)$ is the first Chern class of the line bundle
$T(\cG)$, and the map $T$ on the right hand side is the transfer in integral
cohomology $T:H^*(B,\Z)\rightarrow H^{*-1}(LB,\Z)$.
Modulo torsion  Lemma \ref{com1} follows from the Lemma \ref{com2} below.

If $\cG\in \Gerbe(B)$, then there is a well defined notion of the curvature
$R^{\cG}\in\Omega^3(B)$. This form is closed and represents the image of
$v(\cG)$ in cohomology with real coefficients.

Note that we can define a transfer map $T:\Omega^*(B)\rightarrow
\Omega^{*-1}(LB)$ inducing the transfer in real cohomology
(via the de Rham isomorphism).
 If $\omega\in \Omega^*(B)$,
then we define 
$$T\omega:=\int_{S^1} i_{\partial t} ev^*(\omega) dt\ .$$
\begin{lem}\label{com2}
We have
$2\pi\imath T(R^{\cG})=R^{T(\cG)}$.
\end{lem}
On the right hand side of this equality $R^\cL$ denotes the curvature of the
hermitean line bundle with connection $\cL$. We prove this Lemma at the end of
Subsection \ref{cohomology}.

\subsection{Transgression of the index gerbe}\label{transindex}

Let $\cE$ be a geometric family over $B$ with odd  dimensional fibres.
Then we have an index gerbe  (which was introduced by Lott \cite{lott01})
$\gerbe(\cE)\in\Gerbe(B)$ (we recall the construction in Subsection
\ref{constr}). 

In order to state the main result about the transgression of the index gerbe
we must define the adiabatic limit of the determinant line bundle $\det(
L_\epsilon \cE)$ as $\epsilon\rightarrow 0$.
For $\gamma\in LB$ and $\epsilon\in (0,1]$ let
$p(\gamma,\epsilon):[\epsilon,1]\rightarrow LB\times\R_+$ be the path
$p(\gamma,\epsilon)(t)=(\gamma,t)$.
If $\cL$ is any line bundle with connection over some base and $\sigma$ is a
path in the base defined on $[a,b]\subset \R$, then by $\|^{\cL}_\sigma\in
\End(L_{\sigma(a)},L_{\sigma(b)})$ we denote the parallel transport along
$\sigma$.
The family of maps $\|^{\det(
L \cE)}_{p(\gamma,\epsilon)}$, $\gamma\in LB$, provides an unitary isomorphism
of line bundles $ \Phi_{\epsilon} :\det(
L_\epsilon \cE)\stackrel{\sim}{\rightarrow}  \det(
L_1 \cE)$. It is not compatible with the connections.
We define a family of connections $\nabla^\epsilon:=\Phi_\epsilon\circ
\nabla^{ \det(
L_\epsilon \cE)} \circ \Phi_{\epsilon}^{-1}$ on $\det (L_1\cE)$.
\begin{lem}\label{adiadet}
The limit $\lim_{\epsilon\to 0} \nabla^\epsilon=:\nabla^0$ exists.
\end{lem}
We prove this Lemma in Subsection \ref{adia}.

We define $\det(L_0\cE):=\lim_{\epsilon\to 0} \det(
L_\epsilon \cE):=(\det(
L_1 \cE),h^{\det(L_1\cE)},\nabla^0)$.

We can now state the main result of the present paper.

\begin{theorem}
We have $T \gerbe(\cE) =  \det(
L_0 \cE)$.
\end{theorem}

Since a line bundle with connection is determined by its transgression the
theorem immediately follows from
 \begin{prop}\label{stmain}
We have $T^2 \gerbe(\cE)=T \det(L_0 \cE)$.
\end{prop}

As a first approximation  in Subsection \ref{ccuu} we show the identity of
curvatures \begin{lem}\label{curv} 
$ R^{T\gerbe(\cE)}=R^{ \det(L_0 \cE)}$. 
\end{lem}

The proof of Proposition \ref{stmain} follows from the study of adiabatic
limits of $\eta$-invariants and constitutes the main innovation of the present
paper.

\section{Transgression of gerbes}

\subsection{Gerbes with connection}\label{gerbe}

We describe gerbes with connection on a smooth manifold $B$ in the \v{C}ech
picture (see Hitchin \cite{hitchin99} and Lott \cite{lott01}).

Let $\cU:=(U_\alpha)_{\alpha\in I}$ be an open covering of $B$.  A gerbe $G$ 
on $B$ is given by the following data
$G=((L_{\alpha\beta}),(\theta_{\alpha\beta\gamma}))$:
 \begin{enumerate}
\item A hermitean line bundle $L_{\alpha\beta}$ on each nonempty intersection
$U_\alpha\cap U_\beta$ such that $L_{\beta\alpha}\cong L_{\alpha\beta}^*$.
\item A nowhere vanishing section $\theta_{\alpha\beta\gamma}$ of
$L_{\alpha\beta}\otimes L_{\beta\gamma}\otimes L_{\gamma\alpha}$ over each
nonempty triple intersection $U_\alpha\cap U_\beta\cap U_\gamma$
such that
$\theta_{\beta\gamma\delta} \theta_{\alpha\gamma\delta}^{-1} 
\theta_{\alpha\beta\delta}\theta_{\alpha\beta\gamma}^{-1}=1$ over each
nonempty intersection $U_\alpha\cap U_\beta\cap U_\gamma\cap U_\delta$
(note that this product is a section of a bundle which is canonically trivial).
\end{enumerate}

If $\cV=(V_\mu)_{\mu\in J}$, $s:J\rightarrow I$ is a refinement of $\cU$, then
the same gerbe $G$ is defined by the data
$L^\prime_{\mu\nu}:=(L_{s(\mu)s(\nu)})_{|V_{\mu}\cap V_\nu}$ and
$\theta^\prime_{\mu\nu\sigma}:=(\theta_{s(\mu)s(\nu)s(\sigma)})_{|V_\mu\cap
V_\nu\cap V_\sigma}$.

Two sets of data $G=((L_{\alpha\beta}),(\theta_{\alpha\beta\gamma}))$ and
$G^\prime=((L^\prime_{\alpha\beta}),(\theta^\prime_{\alpha\beta\gamma}))$
describe isomorphic gerbes if there is a family of hermitean line bundles
$L_\alpha$ over $U_\alpha$, $\alpha\in I$ and isomorphisms 
$L^\prime_{\alpha\beta}\cong L_\alpha^{-1} \otimes L_{\alpha\beta}\otimes
L_\beta$ such that $\theta^\prime_{\alpha\beta\gamma}$ corresponds to
$\theta_{\alpha\beta\gamma}$ under the induced isomorphism
$L^\prime_{\alpha\beta}\otimes L^\prime_{\beta\gamma}\otimes
L^\prime_{\gamma\alpha}=L_{\alpha\beta}\otimes L_{\beta\gamma}\otimes
L_{\gamma\alpha}$.
In general, if two gerbes $G,G^\prime$ are defined with respect to coverings
$\cU$, $\cU^\prime$, then they are isomorphic if they are so when defined on
a suitable common refinement. 

A connection of a gerbe $G=((L_{\alpha\beta}),(\theta_{\alpha\beta\gamma}))$
defined with respect to $\cU$ consists of the following data
$((\nabla^{L_{\alpha\beta}}),(F_\alpha))$:
\begin{enumerate}
\item $\nabla^{L_{\alpha\beta}}$ is a hermitean connection on
$L_{\alpha\beta}$ such that
$\nabla^{L_{\beta\alpha}}=(\nabla^{L_{\alpha\beta}})^*$ and such that
$\theta_{\alpha\beta\gamma}$ is parallel with respect to the induced connection
on $L_{\alpha\beta}\otimes L_{\beta\gamma}\otimes
L_{\gamma\alpha}$.
\item
$F_\alpha\in\Omega^2(B)$ satisfies
$F_\beta-F_\alpha=c_1(\nabla^{L_{\alpha\beta}})$ over $U_\alpha\cap U_\beta$,
where $c_1(\nabla):=\frac{-1}{2\pi\imath}
R^{\nabla}$ is the  first Chern form.
\end{enumerate}

If 
$\cV=(V_\mu)_{\mu\in J}$, $s:J\rightarrow I$
is a refinement of $\cU$, then
$\nabla^{L^\prime_{\mu\nu}}:=\nabla^{L_{s(\mu)s(\nu)}}_{ |V_\mu\cap
V_\nu}$ and $F^\prime_\mu:=(F_{s(\mu)})_{|V_\mu}$ define the same connection on
$G$. Two gerbes with connection
$\cG=((L_{\alpha\beta}),(\theta_{\alpha\beta\gamma}),(\nabla^{L_{\alpha\beta}}),(F_\alpha))$
and
$\cG^\prime=((L^\prime_{\alpha\beta}),(\theta^\prime_{\alpha\beta\gamma}),(\nabla^{L^\prime_{\alpha\beta}}),(F^\prime_\alpha))$
are isomorphic, if there is family of hermitean line bundles with connection 
$(L_\alpha,\nabla^{L_\alpha})$ over $U_\alpha$, $\alpha\in I$, such that
$\nabla^{L^\prime_{\alpha\beta}}$ corresponds to the induced
connection $\nabla^{L_\alpha^{-1} \otimes L_{\alpha\beta}\otimes L_\beta}$
and $F^\prime_\alpha=F_\alpha+c_1(\nabla^{L_\alpha})$.
 In general, if two gerbes with connection $\cG,\cG^\prime$ are defined
with respect to coverings $\cU$, $\cU^\prime$, then they are isomorphic if
they are so when defined on a suitable common refinement.

Let
$\cG=((L_{\alpha\beta}),(\theta_{\alpha\beta\gamma}),(\nabla^{L_{\alpha\beta}}),(F_\alpha))$
be a gerbe with connection. 
Then $dF_\alpha=dF_\beta$ over $U_\alpha\cap U_\beta$. Therefore, there is a
form $R^{\cG}\in \Omega^3(B)$ which restricts to $dF_\alpha$ over $U_\alpha$.
This form is called the curvature of $\cG$.  The cohomology class
$[R^{\cG}]\in H^3_{dR}(B)$  is the de Rham cohomology class
corresponding to $v(\cG)\in H^3(B,\Z)$. In particular, it is independent of the
connection.

Let $\Gerbe(B)$ be the set of all isomorphism classes of gerbes with
connection  on $B$. We define a product on $\Gerbe(B)$ as follows.
If $\cG,\cG^\prime \in\Gerbe(B)$ are defined with respect to a covering $\cU$,
then $\cG\otimes \cG^\prime$ is given by
the data
$L_{\alpha\beta}\otimes L^\prime_{\alpha\beta}$,
$\theta_{\alpha\beta\gamma}\otimes\theta^\prime_{\alpha\beta\gamma}$,
the induced connections $\nabla^{L_{\alpha\beta}\otimes
L^\prime_{\alpha\beta}}$ and $F_\alpha+F^\prime_\alpha$.
The trivial element is the gerbe  where  $L_{\alpha\beta}$ is the trivial
bundle, $\theta_{\alpha\beta\gamma}\equiv 1$, $\nabla^{L_{\alpha\beta}}$ is
the trivial connection, and $F_\alpha =0$.
The inverse of $\cG$ is given by $L^{-1}_{\alpha \beta}$,
$\theta_{\alpha\beta\gamma}^{-1}$, the induced connection
$\nabla^{L^{-1}_{\alpha\beta}}$, and $-F_\alpha$.

If $f:B^\prime\rightarrow B$ is a smooth map, then it induces a group
homomorphism $f^*:\Gerbe(B)\rightarrow \Gerbe(B^\prime)$.
If
$\cG=((L_{\alpha\beta}),(\theta_{\alpha\beta\gamma}),(\nabla^{L_{\alpha\beta}}),(F_\alpha)) \in\Gerbe(B)$ is defined with respect to a covering $\cU$,
then
$f^*\cG:=((f^*L_{\alpha\beta}),(f^*\theta_{\alpha\beta\gamma}),(\nabla^{f^*L_{\alpha\beta}}),(f^*F_\alpha))$
is defined with respect to $f^*\cU=(f^{-1}U_\alpha)_{\alpha\in I}$.

\subsection{Transgression}\label{trans}

In this subsection we describe the transgression map
$T:\Gerbe(B)\rightarrow \Line(LB)$.

Let $\cG\in\Gerbe(B)$ be defined with respect to $\cU=(U_\alpha)_{\alpha\in
I}$. We consider on $\Z_n:=\Z/n\Z$ and $S^1$ a cyclic ordering. Let
$t:\Z_n\rightarrow S^1$ be a monotone map and $s:\Z_n\rightarrow I$ be any map.
Then we consider the open set
$V(t,s)\subset LB$ consisting of all loops $\gamma$ such that
$\gamma([t(i),t(i+1)])\in U_{s(i-1)}\cap U_{s(i)}$. The family
$\cV:=(V(t,s))_{(n\in \nat,t:\Z_n\rightarrow S^1, s:\Z_n\rightarrow I)}$ forms
an open covering of $LB$.
We define $T\cG$ by describing the restrictions $T\cG_{|V(t,s)}$ and the
transition maps.

We consider the map $E(t,s):V(t,s)\rightarrow \prod_{i\in \Z_n} U_{s(i-1)}\cap
U_{s(i)}$ given by $E(t,s)(\gamma)=\prod_{i\in\Z_n} \gamma(t(i))$.
Let $p_j:\prod_{i\in \Z_n} U_{s(i-1)}\cap
U_{s(i)}\rightarrow U_{s(j-1)}\cap
U_{s(j)}$ be the canonical projection.
We define the hermitean line bundle
$$T\cG_{|V(t,s)}:=E(t,s)^*\left(\bigotimes_{i\in\Z_n} p_i^*
L_{s(i-1)s(i)}\right)\ .$$

The transition maps are generated by two types.
Assume first that $(t^\prime,s^\prime)$ is a refinement of $(t,s)$. Then there
is a monotone map $r:\Z_n\rightarrow \Z_m$, $m\ge n$ such that
$t=t^\prime\circ r$ and $s^\prime(i)=s(j)$ for all $i\in [r(j),r(j+1))$.
In this case $V(t,s)\subset V(t^\prime,s^\prime)$.
Using the fact that $L_{\alpha\alpha}$ is the trivial bundle we obtain an
isomorphism
$$E(t^\prime,s^\prime)^*\left(\bigotimes_{i\in\Z_m} p_i^*
L_{s^\prime(i-1)s^\prime(i)}\right)_{|V(t,s)}\stackrel{\sim}{\rightarrow} 
 E(t,s)^*\left(\bigotimes_{i\in\Z_n} p_i^* L_{s(i-1)s(i)}\right)\
.$$

The other type of transition map is related to the change of the indexing map
$s$. We fix $t:\Z_n\rightarrow S^1$ and consider two maps
$s,s^\prime:\Z_n\rightarrow I$ such that $V(t,s)\cap
V(t,s^\prime)\not=\emptyset$.
Let $\gamma\in V(t,s)\cap
V(t,s^\prime)$. We want to construct an isomorphism of fibres
$$E(t,s)^*\left(\bigotimes_{i\in\Z_n} p_i^* L_{s(i-1)s(i)}\right)(\gamma)\cong 
E(t,s^\prime)^*\left(\bigotimes_{i\in\Z_n} p_i^*
L_{s^\prime(i-1)s^\prime(i)}\right)(\gamma)\ . $$ 
Using the family of sections
$(\theta_{\alpha\beta\gamma})$ we have isomorphisms \begin{eqnarray*}
\lefteqn{E(t,s)^*\left(\bigotimes_{i\in\Z_n} p_i^*
L_{s(i-1)s(i)}\right)(\gamma)}&&\\ &=&
\bigotimes_{i\in\Z_n} L_{s(i-1)s(i)}(\gamma(t(i)))\\
&=&
\bigotimes_{i\in\Z_n}L_{s(i-1)s(i)}(\gamma(t(i)))\\
&&
\otimes
\bigotimes_{i\in\Z_n} L_{s(i-1)s^\prime(i-1)}(\gamma(t(i)))\otimes
L_{s(i)s(i-1)}(\gamma(t(i)))\otimes L_{s^\prime(i-1)s(i)}(\gamma(t(i)))\\
&&\otimes 
\bigotimes_{i\in\Z_n} L_{s^\prime(i)s(i)}(\gamma(t(i)))\otimes
L_{s^\prime(i-1)s^\prime(i)}(\gamma(t(i)))\otimes
L_{s(i)s^\prime(i-1)}(\gamma(t(i)))\\
&=&\bigotimes_{i\in\Z_n}L_{s^\prime(i-1)s^\prime(i)}(\gamma(t(i)))
\otimes \bigotimes_{i\in\Z_n}
\Hom\left(L_{s(i)s^\prime(i)}(\gamma(t(i))),L_{s(i)s^\prime(i)}(\gamma(t(i+1)))\right)
\end{eqnarray*}
Let $\gamma_i:[t(i),t(i+1)]\rightarrow U_{s(i)}\cap U_{s^\prime(i)}$ be the
path obtained by restriction of $\gamma$.
Then we have the element
$\|^{L_{s(i)s^\prime(i)}}_{\gamma_i}\in
\Hom\left(L_{s(i)s^\prime(i)}(\gamma(t(i))),L_{s(i)s^\prime(i)}(\gamma(t(i+1)))\right)$.
Therefore, $\bigotimes_{i\in\Z_n}\|^{L_{s(i)s^\prime(i)}}_{\gamma_i}$ induces
an ismorphism
\begin{eqnarray*}&&
\bigotimes_{i\in\Z_n}L_{s^\prime(i-1)s^\prime(i)}(\gamma(t(i)))
\otimes \bigotimes_{i\in\Z_n}
\Hom\left(L_{s(i)s^\prime(i)}(\gamma(t(i))),L_{s(i)s^\prime(i)}(\gamma(t(i+1)))\right)\\&
\cong &
\bigotimes_{i\in\Z_n}L_{s^\prime(i-1)s^\prime(i)}(\gamma(t(i)))\ .
\end{eqnarray*}
If we vary $\gamma\in V(t,s)\cap
V(t,s^\prime)$, then we obtain an isomorphism
\begin{eqnarray*}
&&\hspace{-3cm}\Phi((t,s^\prime),(t,s)):E(t,s)^*\left(\bigotimes_{i\in\Z_n}
p_i^* L_{s(i-1)s(i)}\right)_{|V(t,s)\cap
V(t,s^\prime)}\\&\stackrel{\sim}{\rightarrow}&
E(t,s^\prime)^*\left(\bigotimes_{i\in\Z_n} p_i^*
L_{s^\prime(i-1)s^\prime(i)}\right)_{|V(t,s)\cap V(t,s^\prime)}\
.\end{eqnarray*}

We now describe  the connection on the hermitean line bundle constructed above.
Note that the restrictions $T\cG_{|V(t,s)}$ come with induced hermitean
connections.  However, these connections are not compatible with the transition
maps $\Phi((t,s^\prime),(t,s))$. The point is that 
$\gamma\mapsto \bigotimes_{i\in\Z_n}\|_{\gamma_i}^{L_{s(i)s^\prime(i)}}$ is not
parallel. We will fix this problem by introducing a correction using the forms
$F_\alpha$. For our purpose it turns out to be useful to describe the parallel
transport of the connection $\nabla^{T\cG}$. Let $\Gamma:[a,b]\rightarrow
V(t,s)$ be a path. For $i\in\Z_n$ we define the path $\Gamma_i:[a,b]\rightarrow
U_{s(i-1)}\cap U_{s(i)}$ by $\Gamma_i(x)=\Gamma(x)(t(i))$.
Then we define
$$\|^{T\cG}_\Gamma:=\bigotimes_{i\in\Z_n} \|^{L_{s(i-1)s(i)}}_{\Gamma_i}
\prod_{i\in\Z_n} \exp\left(2\pi\imath \int_{[a,b]\times[t(i),t(i+1)]}\Gamma^*
F_{s(i)}\right)\ .$$

We leave it as an exercise to the interested reader to check that the
transgression is well defined. In fact, in Subsection \ref{cohomology} we will
show that this construction coincides with the explicite
description of thetransgression defined in the Deligne cohomology picture.
This also implies well-definedness (see  Corollary \ref{well}).

If $f:B^\prime \rightarrow B$ is a smooth map, then in induces a map
$Lf:LB^\prime\rightarrow LB$.
One can check that transgression is compatible with pull-back, i.e.
\begin{equation}\label{comoa}T\circ f^*=Lf^*\circ T\ ,\end{equation} where
$Lf^*:\Line(LB)\rightarrow \Line(LB^\prime)$ is the pull-back of hermitean
line bundles with connection.

\subsection{The Deligne cohomology picture}\label{cohomology}

If $(R^*,d)$ is a complex of sheaves on a space $B$, then under certain
conditions we can compute its hypercohomology using \v{C}ech cohomology.
If $\cU$ is a covering of $B$, then we form the double complex
$C^{p,q}:=C_\cU^p(R^q)$, $d_0:=\delta$, $d_1:=d$, where
$C_\cU^p(R^q)$ is the space of \v{C}ech $p$-cochains of $R^q$ and $\delta$ is
the differential of the \v{C}ech complex.  
Its action on a \v{C}ech $p$-cochain $X$ is given by
$$\delta X_{i_0\dots i_{p+1}}:=\sum_{j=0}^{p+1} (-1)^j X_{i_0\dots
\check{i}_j\dots i_{p+1}}$$
in the usual notation. We write this down in order to fix our sign conventions.
The hypercohomology of the complex $(R^*,d)$ is approximated (one has to go
to the limit over all coverings) by the total complex complex
$(C_{tot}^*,d_{tot})$ associated to the double complex above.
In order to fix signs, $d_{tot} c^{p,q}:= \delta c^{p,q}+ (-1)^{p} d c^{p,q}$
for $c^{p,q}\in C^{p,q}$.

Let $\cK^p_B$ be the complex of sheaves on a manifold $B$
$$0\rightarrow
\underline{U(1)}\stackrel{\frac{1}{2\pi\imath}
d\log}{\rightarrow}\Omega^1\stackrel{d}{\rightarrow}\dots
\stackrel{d}{\rightarrow}\Omega^p\rightarrow 0\ ,$$ where $\underline{U(1)}$
denotes the sheaf of functions with values in $U(1)$, and $\Omega^q$ is the
sheaf of real $q$-forms. By $H^k(B,\cK^p_B)$ we denote its $k$'th
hypercohomology. There is an integration over the fibre  $$\int_{S^1}
H^k(B\times S^1,\cK^p_{B\times S^1})\rightarrow H^{k-1}(B,\cK^{p-1}_B)$$ which
was introduced by Gawedzki and is described in the \v{C}ech picture in
Brylinski`s book \cite{brylinski93}, Proposition  6.5.2. We can now define the
transgression $$T:H^k(B,\cK^p_B)\rightarrow H^{k-1}(LB,\cK^{p-1}_{LB})$$
as the composition $T:=\int_{S^1}\circ ev^*$.

There are natural isomorphisms
$u_B:C^\infty(B,U(1)) \stackrel{\sim}{\rightarrow} H^0(\cK^0_B)$,
$l_B:\Line(B)\stackrel{\sim}{\rightarrow} H^1(B,\cK^1_B)$ and
$g_B:\Gerbe(B)\stackrel{\sim}{\rightarrow}  H^2(B,\cK^2_B)$. We will show that
the transgression in Deligne cohomology is compatible with the transgression
defined on the level of geometric objects.

First we describe the maps $u_B,l_B,g_B$ explicitly. 
A function $f\in C^\infty(B,U(1))$ represents a cohomolohy class $u_B(f)\in
H^0(B,\cK^0_B)=H^0(B,\underline{U(1)})$ in a natural way.

Let $\cL=(L,h^L,\nabla^L)\in\Line(B)$. We choose a covering
$\cU:=(U_\alpha)_\alpha\in I$ such that $L_{|U_\alpha}$ is trivial.
We fix sections $s_\alpha\in C^\infty(U_\alpha,L)$ of unit length for all
$\alpha\in I$. Then we define $A_\alpha:=\nabla^{L}\log s_\alpha$,
$A_\alpha\in\imath \Omega^1(U_\alpha)$.
Let $U_{\alpha\beta}\in C^\infty(U_\alpha\cap U_\beta, U(1))$
be the transition function $U_{\alpha\beta}:=\frac{s_\alpha}{s_\beta}$. Then
$((U_{\alpha\beta}),(\frac{-1}{2\pi\imath}A_\alpha))$ is a \v{C}ech cocylcle
representing the  hypercohomology class $l_B(\cL)$. Note that over $U_\alpha\cap U_\beta$ 
we have $A_\alpha-A_\beta=d\log
U_{\alpha\beta}$.

Let $\cG=((L_{\alpha\beta}),(\theta_{\alpha\beta\gamma}),(\nabla^{L_{\alpha\beta}}),(F_\alpha))$
be a gerbe with connection on $B$ defined with respect to $\cU$.
We assume that $L_{\alpha\beta}$ is trivial, and we fix unit length sections
$s_{\alpha\beta}$. Then we define
$A_{\alpha\beta}:=\nabla^{L_{\alpha\beta}}\log s_{\alpha\beta}$,
$A_{\alpha\beta}\in\imath \Omega^1(U_\alpha\cap U_\beta)$. Furthermore, we can
identify the section $\theta_{\alpha\beta\gamma}$ with a function
$\frac{\theta_{\alpha\beta\gamma}}{s_{\alpha\beta}s_{\beta\gamma}s_{\gamma\alpha}}\in
C^\infty(U_\alpha\cap U_\beta \cap U_\gamma,U(1))$. Then
$((\theta_{\alpha\beta\gamma}),(\frac{1}{2\pi\imath}
A_{\alpha\beta}),(-F_\alpha))$ is  a  \v{C}ech cocycle
representing the  hyper cohomology class $g_B(\cG)$.

\begin{lem}
$T\circ l_B= u_{LB}\circ T$.
\end{lem}
\proof
Let $\cL=(L,h^L,\nabla^L)\in\Line(B)$.
Then $u_{LB}\circ T(\cL)$ is just the $U(1)$-valued function which maps
$\gamma\in LB$ to the holonomy $\hol_\gamma(\cL)$ of $\cL$ along $\gamma$.
Let $l_B(\cL)$ be
represented by the \v{C}ech cocycle
$((U_{\alpha\beta}),(\frac{-1}{2\pi\imath}A_\alpha))$ as above.
We consider $\gamma\in LB$. We choose $n\in\nat$, a monotone map
$t:\Z_n\rightarrow S^1$ and $s:\Z_n\rightarrow I$ such that
$\gamma([t(i),t(i+1)])\subset U_{s(i)}$. We choose small neighbourhoods $W_i$
of $[t(i),t(i+1)]$ such that $\gamma(\bar W_i)\subset U_{s(i)}$.
We define the neighbourhood $\tilde V$ of $\gamma$ to be the space
of all $\gamma^\prime\in LB$ satisfying $\gamma^\prime(\bar W_i)\subset
U_{s(i)}$. Then $ev(\tilde V\times W_i)\subset U_{s(i)}$ for all $i\in \Z_n$.
The restriction of $ev^* l_B(\cL)$
is represented by $((u_{ij}),(\frac{-1}{2\pi\imath} a_i))$ where $u_{ij}$ is
defined if $j=i\pm 1$  by $u_{ (i+1)i}(\gamma,x)=U_{s(i+1)s(i)}(\gamma(x))$,
$x\in W_i\cap W_{i+1}$, and $a_i:=  ev^* A_{s(i)} \in
\Omega^1(\tilde V\times W_i)$.  
The description of $\int_{S^1}$ given in \cite{brylinski93}, Equation (6.21), 
yields 
$$T\circ l_B(\cL)(\gamma)=\prod_{i\in \Z_n} \exp\left(-
\int_{[t(i),t(i+1)]} i_{\partial t} a_i dt\right) \prod_{j\in \Z_n}
u_{j,j+1}(\gamma,t(i+1))^{-1}\ .$$
The right hand side can be rewritten as
$$\prod_{i\in \Z_n} \exp\left(-
\int_{[t(i),t(i+1)]} i_{\partial \gamma_t} A_{s(i)} dt\right) \prod_{j\in \Z_n}
U_{s(j+1)s(j)}(\gamma(t(i+1)))
=\hol_\gamma(\cL)\ .$$
\hB

\begin{lem}
$T\circ g_B =  l_{LB}\circ T$
\end{lem}
\proof
It suffices to show that
$T^2 \circ g_B= T\circ l_{LB}\circ T$
holds true in $H^0(L^2B,\cK^0_{L^2B})\cong C^\infty(L^2B,U(1))$.
Let 
$\cG=((L_{\alpha\beta}),(\theta_{\alpha\beta\gamma}),(\nabla^{L_{\alpha\beta}}),(F_\alpha))$
be a gerbe with connection on $B$ defined with respect to $\cU$.
In Subsection \ref{trans} we described the parallel transport of $T\cG$.
Let $\Gamma\in L^2B$. Then we choose a monotone map $u:\Z_m\rightarrow S^1$
such that there is a fixed monotone function $t:\Z_n\rightarrow S^1$ and
a family $s_i:\Z_n\rightarrow I$ such that
$\Gamma([u(i),u(i+1)])\subset V(t,s_i)$ for all $i\in\Z_n$. 
We define $\gamma_{ij}:[u(i),u(i+1)]\rightarrow B$ by
$\gamma_{ij}(x)=\Gamma(x,t(j))$.
Then we have 
\begin{eqnarray*}
\hol_{\Gamma}(T\cG)&=&{\Huge \circ}_{i\in \Z_n}
\|^{T\cG}_{\Gamma_{|[u(i),u(i+1)]}}\\&=&
{\Huge \circ}_{i\in \Z_n} \left(\Phi((t,s_{i+1}),(t,s_i))\circ 
\bigotimes_{j\in\Z_m} \|^{L_{s_i(j-1)s_i(j)}}_{\gamma_{ij}} \right)\\&&
\prod_{i\in \Z_n}\prod_{j\in\Z_m} \exp\left(2\pi\imath
\int_{[u(i),u(i+1)]\times[t(j),t(j+1)]}\Gamma^* F_{s_i(j)}\right)
\end{eqnarray*}
Writing out the first factor and the transition maps $\Phi$ with respect to the
given trivializations of the $L_{\alpha\beta}$  we get
\begin{eqnarray*}
\hol_{\Gamma}(T\cG)&=&
\prod_{i\in \Z_n} \prod_{j\in\Z_m}
\theta_{s_i(j-1)s_{i+1}(j-1)s_i(j)}(\Gamma(u(i+1),t(j)))
\theta_{s_i(j)s_{i+1}(j-1)s_{i+1}(j)}(\Gamma(u(i+1),t(j))) \\&&
\prod_{i\in \Z_n} \prod_{j\in\Z_m}\exp\left(- \int_{[t(j),t(j+1)]}
\Gamma(u(i+1))^* A_{s_i(j)s_{i+1}(j)}\right)\\ &&\prod_{i\in
\Z_n}\prod_{j\in\Z_m} \exp \left(- \int_{[u(i),u(i+1)]} \Gamma(.,t(j))^*
A_{s_i(j-1)s_i(j)}\right) \\&&\prod_{i\in \Z_n}\prod_{j\in\Z_m}
\exp\left(2\pi\imath \int_{[u(i),u(i+1)]\times[t(j),t(j-1)]}\Gamma^*
F_{s_i(j)}\right) \end{eqnarray*}

We now describe $T\circ g_B(\cG)$.
We choose neighbourhoods $W_j$ of $[t(j),t(j+1)]$ and
$R_i$ of $[u(i),u(i+1)]$ such that $\Gamma(\bar R_i\times \bar W_j)\subset
U_{s_i(j-1)}\cap U_{s_i(j)}$. Let $\tilde V_i\subset LB$ be the set of all
$\gamma$ such that $\gamma(\bar W_j)\subset U_{s_i(j-1)}\cap U_{s_i(j)}$
for all $j\in\Z_m$.
We now describe the restriction of $T\circ g_B(\cG)$ to $\cup_{i\in
\Z_n} \tilde V_i$.
The restriction of $ev^* g_B(\cG)$ to $\cup_{i\in
\Z_n,j\in \Z_m} \tilde V_i \times W_j$ is given by
$$((ev^* \theta_{s_{i_1}(j_1)s_{i_2}(j_2)s_{i_3}(j_3)}
),(\frac{1}{2\pi\imath}ev^* A_{s_{i_1}(j_1)s_{i_2}(j_2)}),(-ev^*F_{s_i(j)}))
\ . $$
Applying $\int_{S^1}$ we obtain a cocycle
$((X_{i_1i_2}),(Y_{i}))$ representing the restriction of $T\circ g_B(\cG)$,
where \begin{eqnarray*}
X_{i_1i_2}(\gamma) &=&
\prod_{j\in\Z_m}\theta_{s_{i_1}(j-1)s_{i_1}(j)s_{i_2}(j)}(\gamma(t(j))
\theta_{s_{i_1}(j-1)s_{i_2}(j-1)s_{i_2}(j)}(\gamma(t(j))^{-1}\\
&&\prod_{j\in\Z_m}
\exp\left(\int_{[t(j),t(j+1)]}\gamma^*A_{s_{i_1}(j)s_{i_2}(j)})  \right)\\
Y_{i}&=&-\sum_{j\in\Z_m}  \int_{[t(j),t(j+1)]} i_{\partial_t}ev^* F_{s_i(j)}
dt\\
&&-\sum_{j\in \Z_m} ev_{|.\times \{t(j+1)\}}^* A_{s_{i}(j)s_{i}(j+1)}\ .
\end{eqnarray*}

It now suffices to compute $ev^*_{\{\Gamma\}\times S^1}((X_{i_1i_2}),(Y_{i}))$
with respect to the covering $(\{\Gamma\}\times \tilde V_i)_{i\in\Z_n}$ and
apply $\int_{S^1}$ again. The result is then
\begin{eqnarray*}
\lefteqn{T^2\circ g_B(\cG)}&&\\
&=&
\prod_{i\in\Z_n}\prod_{j\in\Z_m}\theta_{s_i(j-1)s_i(j)s_{i+1}(j)}
(\Gamma(u(i+1),t(j))^{-1}
\theta_{s_i(j-1)s_{i+1}(j-1)s_{i+1}(j)}(\Gamma(u(i+1),t(j))\\
&&\prod_{i\in\Z_n}\prod_{j\in\Z_m} \exp\left(  - \int_{[t(j),t(j+1)]}
\Gamma(u(i+1))^* A_{s_i(j)s_{i+1}(j)} \right) \\&&
\prod_{i\in\Z_n}\prod_{j\in\Z_m}   \exp \left( -\int_{[u(i),u(i+1)]}
\Gamma(.,t(j+1))^* A_{s_i(j-1)s_i(j)}\right)\\&&
\prod_{i\in\Z_n}\prod_{j\in\Z_m}\exp\left(2\pi\imath
\int_{[u(i),u(i+1)]\times [t(j),t(j+1)]}\Gamma^*F_{s_i(j)} \right)
\end{eqnarray*}
Using that $\delta (\theta_{\alpha\beta\gamma})=0$ it is now easy to see that
$T^2\circ g_B(\cG)(\Gamma)=\hol_{\Gamma}(T\cG)$.
\hB

An immediate consequence of this proof is
\begin{kor}\label{well}
The transgression $T:\Gerbe(B)\rightarrow \Line(L)$
is well-defined by the construction given in Subsection \ref{trans}.
\end{kor}

It now follows from \cite{brylinski93}, Proposition  6.5.1, that the transgression of
gerbes with connection constructed in \ref{trans} is equivalent to the
negative of the construction given in \cite{brylinski93}, 6.2.1, which is due
to Deligne and Brylinski.

The equation \cite{brylinski93}, (6.8) for the curvature of 
$T\cG$ shows
$$R^{T\cG}=2\pi\imath T R^{\cG}\ .$$
This proves the assertion of Lemma \ref{com2}. \hB

\section{The determinant line bundle}

\subsection{Adiabatic limit of determinant line bundles}\label{adia}

The goal of this subsection is the proof of Lemma \ref{adiadet}.
We first recall the Bismut-Freed  formula for the curvature of the determinant
line bundle $\det(\cE)$  associated to a geometric family
$\cE$ over some base $B$ with even  dimensional fibres.
\begin{equation}\label{detcurv}R^{\det(\cE)}= \left[\int_{M/B}
\hA(\nabla^{T^v\pi})\ch(\nabla^V)\right]_{(2)}\ . \end{equation} Here
$[.]_{(2)}$ stands for the two form component. For a proof of this formula we
refer to  \cite{berlinegetzlervergne92}, 10.35. Note that in the present paper
the Chern character form $\ch$ and the $\hA$-genus already include the
$2\pi\imath$-factors (we use the topologist's convention) whereas in
\cite{berlinegetzlervergne92} they are not included.

Let now $\cE$ be a geometric family with odd dimensional fibres and consider
the loop $L\cE$ over $LB\times \R_+$.  The form $\hA(\nabla^{T^vL\pi})$
depends polynomially on $\epsilon\in \R_+$. In particular, its limit as
$\epsilon\to 0$ is $\hA(\nabla^{T^v \tilde ev^*\pi})$ using the notation
introduced in Subsection \ref{geo}.

Recall the definition of the connection $\nabla^\epsilon$ on $\det(L_1\cE)$
given in Subsection \ref{transindex}. We can express the parallel transport
$\|^{\nabla^\epsilon}$ of this connection in terms of the parallel transport 
$\|^{\det (L_1\cE)}$ as follows. Let $\Gamma:[a,b]\rightarrow LB$ be any path.
Then
$$\|_\Gamma^{\nabla^\epsilon}=\|_\Gamma^{\det (L_1\cE)}\exp\left(-
\int_{[\epsilon,1]\times [a,b]} \Psi^* R^{\det{L\cE}}\right)\ ,$$
where $\Psi:[\epsilon,1]\times [a,b]\rightarrow LB\times \R_+$ is given by 
$\Psi(\delta,x)=(\Gamma(x),\delta)$.
We now employ (\ref{detcurv}) and perform the limit $\epsilon\to 0$.
We obtain $$\|_\Gamma^{\nabla^0}=\|_\Gamma^{\det (L_1\cE)}\exp\left(-
\int_{(0,1]\times [a,b]} \Psi^*  \left[\int_{ev^* M/B}
\hA(\nabla^{T^v L\pi})\ch(\nabla^{Ev^* V})\right]_{(2)}\right)\ ,$$
where the integral converges in view of the remark above. \hB

\subsection{Adiabatic limit of the holonomy}

In Subsection \ref{adia} we have seen that the limit $\lim_{\epsilon\to
0}\nabla^\epsilon=:\nabla^0$ exists. Recall the definition $\det(L_0\cE)=(\det(
L_1 \cE),h^{\det(L_1\cE)},\nabla^0)$.
By Theorem \ref{bf} we know that $T\det(L_\epsilon\cE)=\lim_{\delta\to 0}
\tau(D(L_\delta L_\epsilon\cE))$ and therefore
$$T \det (L_0\cE)=\lim_{\epsilon\to 0} \lim_{\delta\to 0}
\tau(D(L_\delta L_\epsilon\cE))\ .$$
 \begin{lem}\label{diag}
We have $T \det (L_0\cE)=\lim_{\epsilon\to 0} \tau(D(L_\epsilon
L_\epsilon\cE))$.
\end{lem}
\proof
We fix $\Gamma\in L^2(B)$. Then we define the map
$f:\R_+\times\R_+\rightarrow L(LB\times \R_+)\times\R_+$ by
$f(\epsilon,\delta)(t)=((\Gamma(t),\delta),\epsilon)$.
We consider the geometric family $f^*L^2\cE$ over
$\R_+\times\R_+$.
Then we have $\tau(D(L_\epsilon
L_\delta\cE))(\Gamma)=\tau(D(f^*L^2\cE))(\epsilon,\delta)$.
The differential of the $\tau$-function is given by \cite{daifreed94}, Thm.
1.9,
$$d\tau(D(f^*L^2\cE))=2\pi\imath \left[\int_{fibre}
\hA(\nabla^{T^vf^*L^2\pi})\ch(\nabla^\cW)\right]_{(1)}\ ,$$
where $\int_{fibre}$ means integration over the fibre of the family
$f^*L^2\cE$, and $\cW$ is the induced hermitean bundle with connection of the
family $f^*L^2\cE$.

Again, $\hA(\nabla^{T^vf^*L^2\pi})$ depends polynomially on $\epsilon,\delta$
and is in particular uniformly bounded for small $(\epsilon,\delta)$.
This justifies the change of the path of integration below.
Let $\gamma,\gamma^\prime$ be the path`s in $\R_+\times\R_+$ given by
$\gamma(t):=(1-2t,1)$, $t\in [0,1]$, and $\gamma(t):=(0,1-2(t-\frac12))$ for
$t\in [\frac12,1]$,  $\gamma^\prime(t):=(1-t,1-t)$, $t\in[0,1]$.
Then we can compute
\begin{eqnarray*}
\lim_{\epsilon\to 0} \lim_{\delta\to 0}
\tau(D(L_\delta L_\epsilon\cE))(\Gamma)&=& \tau(D(L_1 L_1\cE))(\Gamma) +
\int_\gamma d\tau(D(f^*L^2\cE))\\
&=& \tau(D(L_1 L_1\cE))(\Gamma)+ \int_{\gamma^\prime}
d\tau(D(f^*L^2\cE))\\
&=&\lim_{\epsilon\to 0} \tau(D(L_\epsilon
L_\epsilon\cE))(\Gamma)\ .
\end{eqnarray*}
\hB

\section{The index gerbe}

\subsection{Construction of the index gerbe}\label{constr}

Let $\cE$ be a geometric family over some base $B$ with odd dimensional
fibres. Furthermore, let $Q\in \End(\Gamma(\cE))$ be a family of
finite dimensional selfadjoint operators such that
$Q_b=E_{D_b(\cE)}(-R(b),R(b))Q_b E_{D_b(\cE)}(-R(b),R(b))$, $b\in B$, for some
continuous function $R\in C(B)$, where $E_{D_b(\cE)}$ denotes the spectral
projection.

If $D(\cE)+Q$ is invertible for every $b\in B$, then following
\cite{daizhang96} we call $Q$  a Melrose-Piazza operator. 
It was shown by Melrose-Piazza in \cite{melrosepiazza97} that a Melrose-Piazza
operator exists iff $\ind(\cE)=0$.

The index gerbe $\gerbe(\cE)\in\Gerbe(B)$ was introduced by Lott
\cite{lott01}.  Here we recall this definition and discuss a slight
generalization $\gerbe(\cE,Q)$, which is the index gerbe of the perturbed
family $D(\cE,Q):=D(\cE)+Q$.    

Let $\cU=(U_\alpha)_{\alpha\in I}$ be a covering of $B$.
If we choose these sets sufficiently small then there exists functions
$h_\alpha\in C^\infty_c(\R)$, $\alpha\in I$, such that
$D_{b,\alpha}:=D_b(\cE,Q)+h_\alpha(D_b(\cE,Q))$ is invertible for all $b\in
U_\alpha$. In fact, in order to construct $h_\alpha$ it suffices to find a
family of open intervals $I_\alpha , J_\alpha\subset \R_+$, $\alpha\in I$, such
that $(I_\alpha\cup - J_\alpha) \cap
\spec(D_b(\cE,Q))=\emptyset$ for all $b\in U_\alpha$. Then one takes for
$h_\alpha$ a function which is zero outside $[-J_\alpha,I_\alpha]$, and is
equal to one in the region between $ -J_\alpha$ and $I_\alpha$  

With this defnition we can now repeat Lott's definition and arguments line by
line. The only difference is that \cite{lott01}, (2.4), holds if one
replaces "$\lim$" by "$\LIM$". But the conclusion
\cite{lott01}, Proposition 1., is still valid.  

Let us recall the construction of
$\gerbe(\cE,Q)$.
 Let $E_{\alpha\beta}^\pm\rightarrow U_\alpha\cap U_\beta$
denote the $\pm 2$-eigenspace bundle of
$\frac{D_\beta}{|D_\beta|}-\frac{D_\alpha}{|D_\alpha|}$.
This bundle comes with a natural hermitean metric and a connection induced by
$\nabla^{\Gamma(\cE)}$.
We define the hermitean line
bundle $L_{\alpha\beta}:=\det(E_{\alpha\beta}^+)\otimes
\det(E_{\alpha\beta}^-)^{-1}$,  and we let $\nabla^{L_{\alpha\beta}}$ be the
induced connection.
Over $U_\alpha\cap U_\beta\cap U_\gamma$ the product $L_{\alpha\beta}\otimes
L_{\beta\gamma}\otimes L_{\gamma\alpha}$ is canonically trivial. 
This gives the sections $\theta_{\alpha\beta\gamma}$. Finally, we define the
two forms $F_\alpha:=\left[\hat \eta_\alpha\right]_{(2)}$,
where
$$\hat \eta_\alpha=\cR\: \LIM_{t\to 0}\frac{1}{\pi^{1/2}} \int_t^\infty
\Tr_\sigma \left(\frac{dA_s}{ds} \ee^{-A^2_s} \right) ds\ .$$
with the super connection 
$$A_s=s\sigma D_\alpha + \nabla^{\Gamma(\cE)} + \frac{\sigma}{4s} c(T)$$
(see Subsection \ref{dirac} for notation), and
$\cR\in\End(\Omega(B))$ multiplies $p$-forms by $(2\pi\imath)^{-{p/2}}$.

The index gerbe $\gerbe(\cE,Q)$ is now given by 
$((L_{\alpha\beta}),(\theta_{\alpha\beta\gamma}),(\nabla^{L_{\alpha\beta}}),(F_\alpha))$.
Formally, $\gerbe(\cE):=\gerbe(\cE,0)$.

The curvature of $\gerbe(\cE,Q)$ is independent of $Q$ and given by
\cite{lott01}, Thm. 1, 
\begin{equation}\label{curva} R^{\gerbe(\cE,Q)}=\left[\int_{M/B}
\hA(\nabla^{T^v\pi})\ch(\nabla^{V})\right]_{(3)}\in\Omega^3(B)\ .\end{equation}

Note that $\gerbe(\cE,Q)$ may depend on $Q$. But its transgression does not.
\begin{lem}\label{indep}
We have $T\gerbe(\cE,Q)=T\gerbe(\cE)$.
\end{lem}
\proof
Since the transgression of line bundles with connections is injective it
suffices to show that $T^2\gerbe(\cE,Q)=T^2\gerbe(\cE)$.
 
Let $\Gamma\in L^2B$. Then $\Gamma:S^1\times S^1\rightarrow B$. 
Furthermore, let $1\in L^2(S^1\times S^1)$ be given by the identity map. Then
by (\ref{comoa}) 
we can write $T^2\gerbe(\cE,Q)(\Gamma)=T^2\gerbe(\Gamma^*\cE,\Gamma^*Q)(1)$.

Let $\pr:S^1\times S^1\times \R \rightarrow S^1\times S^1$ be the projection.
We define $\tilde Q\in \End(\Gamma(\pr^* \Gamma^*\cE))$ such that $\tilde
Q_{|S^1\times S^1 \times \{\epsilon\}}=\epsilon \Gamma^* Q$. By (\ref{curva})
we have on the one hand $R^{\gerbe(\pr^*\Gamma^*\cE,\tilde
Q)}=\pr^*R^{\gerbe(\Gamma^*\cE)}$. On the other hand,
$R^{\gerbe(\Gamma^*\cE)}=0$ since $\dim(S^1\times S^1)=2<3$.

Furthermore, $H^3(S^1\times S^1\times \R,\Z)=0$.
Therefore $g_{S^1\times S^1\times \R}(\gerbe(\pr^*\Gamma^*\cE,\tilde Q))$ can
be represented by a cocycle of the form $((1),(0),(-F_\alpha))$, where
$F_\alpha$ is the restriction of a global two form which is closed since the
curvature of $\gerbe(\pr^*\Gamma^*\cE,\tilde
Q)$ is trivial. It is now obvious that
\begin{eqnarray*}
 T^2\gerbe(\Gamma^*\cE,\Gamma^*Q)(1)&=&\exp\left(2\pi\imath \int_{S^1\times
S^1\times\{1\}} F\right)\\
&=&\exp\left(2\pi\imath \int_{S^1\times
S^1\times\{0\}} F\right)\\
&=&T^2\gerbe(\Gamma^*\cE)(1)
\end{eqnarray*}
Therefore, $T^2\gerbe(\cE,Q)(\Gamma)=T^2\gerbe(\cE)(\Gamma)$. \hB

It follows immediately from the definitions that if $\cE=\cE_0\cup_B\cE_1$,
then $\gerbe(\cE)=\gerbe(\cE_0)\otimes\gerbe(\cE_1)$.

\subsection{An identity of curvatures}\label{ccuu}

In this Subsection we prove Lemma \ref{curv}.
By equation (\ref{curva}) and Lemma \ref{com2} we have
\begin{eqnarray*}R^{T\gerbe(\cE)}&=&2\pi\imath T R^{\gerbe(\cE)}\\
 &=&
2\pi\imath \int_{LB\times S^1/LB}   ev^*  \left[ \int_{M/B}
\hA(\nabla^{T^v\pi})\ch(\nabla^{V})\right]_{(3)}\\
&=&2\pi\imath\left[\int_{ev^*M/LB} \hA(\nabla^{T^v ev^*\pi})\ch(\nabla^{Ev^*
V})\right]_{(2)}\ . \end{eqnarray*}
The curvature of $\det(L_\epsilon\cE)$ is given by
$$R^{\det(L_\epsilon\cE)}:= I_\epsilon^* \left[\int_{\tilde{ev}^*M/LB\times
\R_+} \hA(\nabla^{T^v\tilde L\pi})\ch(\nabla^{\tilde{Ev}^*V})\right]_{(2)}\ .$$
Since $\lim_{\epsilon\to 0} \hA(\nabla^{T^vL\pi})=\hA(\nabla^{T^v ev^*\pi})$ we
conclude that $R^{\det( L_0\cE)}=\lim_{\epsilon\to 0}R^{\det(L_\epsilon\cE)}=
2\pi\imath T R^{\gerbe(\cE)}$. \hB

\subsection{Holonomy. Bounding case}

Consider again a  loop of loops $\Gamma\in L^2B$, i.e. a map $\Gamma:S^1\times
S^1\rightarrow B$. We say that $\Gamma$ bounds, iff $[\Gamma:S^1\times
S^1\rightarrow B]$ is trivial in the bordism group
$\Omega^{spin}_2(B)$, i.e. there is a $3$-dimensional oriented spin
manifold  $W$ such that $\partial W\cong S^1\times S^1$ as oriented
spin manifolds, and such that $\Gamma$ extends to $\tilde \Gamma:W\rightarrow
B$. 

We now show the following partial case of Proposition \ref{stmain}.
 
\begin{prop}\label{bounding}
If $\Gamma$ bounds, then
$T^2(\gerbe(\cE))(\Gamma)=T(\det(L_0\cE))(\Gamma)$.
\end{prop}
\proof
We employ the extension $\tilde \Gamma:W\rightarrow
B$ in order to compute both sides in
terms of curvatures.

Let $1\in L^2(S^1\times S^1)$ be given by the identity. Then we have
$T^2(\gerbe(\cE))(\Gamma)=T^2(\gerbe(\Gamma^*\cE))(1)$.
Now $\gerbe(\Gamma^*\cE)$ is the restriction of $\gerbe(\tilde\Gamma^*\cE)$
to the boundary of $W$. Since $H^3(W,\Z)=0$ we have
$$T^2(\gerbe(\Gamma^*\cE))(1)=\exp\left(2\pi\imath \int_{W}
R^{\tilde\Gamma^*\cE}\right)\ .$$ By (\ref{curva}) we have
$$T^2(\gerbe(\cE))(\Gamma)=\exp\left(2\pi\imath \int_{W}  \tilde \Gamma^*
\left[\int_{M/B}
\hA(\nabla^{T^v\pi})\ch(\nabla^{V})\right]_{(3)}\right)\ .$$

We now compute $$T(\det(L_0\cE))(\Gamma)=\lim_{\epsilon\to 0}
T(\det(L_\epsilon\cE))(\Gamma)=\lim_{\epsilon\to 0}\lim_{\delta\to 0}
\tau(D(L_\delta L_\epsilon \cE))(\Gamma)\ ,$$
and hence by Lemma \ref{diag},
$T(\det(L_0\cE))(\Gamma)=\lim_{\epsilon\to 0}
\tau(D(L_\epsilon L_\epsilon \cE))(\Gamma)$.
In other words, $T(\det(L_0\cE))(\Gamma)$
is equal to the adiabatic limit $\lim_{\epsilon\to 0}\tau(D_\epsilon)$, where 
$D_\epsilon$ is the total Dirac operator on $\Gamma^* M$
(see Subsection \ref{dirac}). 
We choose $\tilde\Gamma$ such that it is independent of the normal variable
in a tubular neighbourhood of $\partial W$. Furthermore, we choose a
Riemannian metric $g^W$ as a product in that neighbourhood. Then
the total Dirac operator $\tilde D_\epsilon$ on $\tilde \Gamma^*M$ has a
product structure at $\partial \tilde \Gamma^*M$.

The index theorem of Atiyah-Patodi-Singer \cite{atiyahpatodisinger75} gives
$$\tau(D_\epsilon)=\exp\left(2\pi\imath
\int_{\tilde\Gamma^*M}
\hA(\nabla_\epsilon^{T\tilde\Gamma^*M})\ch(\nabla^{\tilde\Gamma^*V})\right)\
,$$ where $\nabla_\epsilon^{T\tilde\Gamma^*M}$ is the Levi Civita connection
to the metric $g^{T\tilde\Gamma^*M}_\epsilon$.
Now $\lim_{\epsilon\to 0}
\hA(\nabla_\epsilon^{T\tilde\Gamma^*M})=\hA(\nabla^{T^v\tilde \Gamma^*\pi})$.
We conclude that
\begin{eqnarray*}T(\det(L_0\cE))(\Gamma)&=&\exp\left(2\pi\imath
\int_{\tilde\Gamma^*M}
\hA(\nabla^{T^v\tilde \Gamma^*\pi})\ch(\nabla^{\tilde\Gamma^*V})\right)\\&=&
\exp\left(2\pi\imath \int_{W}  \tilde \Gamma^*
\left[\int_{M/B}
\hA(\nabla^{T^v\pi})\ch(\nabla^{V})\right]_{(3)}\right)\ .
\end{eqnarray*}
\hB

\subsection{Holonomy. Trivial $K^1$-case.}

In this subsection we consider another special case of Proposition
\ref{stmain}. Again, let $\Gamma\in L^2B$.

\begin{prop}\label{k1t}
If $\Gamma^*\ind(\cE)=0$,
then $T^2(\gerbe(\cE))(\Gamma)=T(\det(L_0\cE))(\Gamma)$. 
\end{prop}
\proof
In the present subsection we will prove this proposition up to a
certain assertion about adiabatic limits of $\eta$-invariants Proposition 
\ref{limit}.

Again, we write $T^2(\gerbe(\cE))(\Gamma)=T^2(\gerbe(\Gamma^*\cE))(1)$.
Since $\ind(\Gamma^*\cE)=0$ we can find a Melrose-Piazza operator
$Q\in\End(\Gamma(\Gamma^*\cE))$ (see Subsection \ref{constr}).
By Lemma \ref{indep} we have
\begin{equation}\label{step1}
T^2(\gerbe(\Gamma^*\cE))(1)=T^2(\gerbe(\Gamma^*\cE,Q))(1)\ .
\end{equation}
Since $D(\Gamma^*\cE,Q)$ is invertible we can represent
$\gerbe(\Gamma^*\cE,Q)$ with respect to the covering
$\cU=(S^1\times S^1)$, i.e. by the global two form $F:=\left[\hat
\eta\right]_{(2)}$, the two form component of the eta form of $  
D(\Gamma^*\cE,Q)$. In particular, we have
\begin{equation}\label{zw2}T^2(\gerbe(\Gamma^*\cE,Q))(1)=\exp\left(2\pi\imath
\int_{S^1\times S^1} F \right) \ .\end{equation} 
Let $D_\epsilon$ be the total Dirac operator on $\Gamma^*\cE$ with respect to
the metric $g^{T\Gamma^*M}_\epsilon$.
Then as in the proof of Proposition  \ref{bounding}
we have 
\begin{equation}\label{zwi}T(\det(L_0\cE))(\Gamma)=\lim_{\epsilon\to
0}\tau(D_\epsilon)\ .\end{equation} 

We choose some function $\rho\in C^\infty(\R_+)$ such that $\rho(s)=0$ for
$s\le 1$ and $\rho(s)=const$ for $t\ge 2$. 
Let $D_\epsilon(s)=sD_\epsilon+s\rho(s) Q$.
We define
$$\eta(\epsilon,\rho):=\frac{2}{\pi^{1/2}}\int_{0}^\infty \Tr
\frac{d}{ds}D_\epsilon(s) \ee^{-D_\epsilon(s)^2} ds\ .$$
Since $\rho(s)\equiv 0$ for small $s$ this integral converges at $s=0$.
Moreover, since $\rho(s)$ is constant for large $s$ the integrand
vanishes exponentially if $s\to\infty$.
We further define
$$\tau(\epsilon,\rho):=\exp\left(2\pi\imath\frac{\eta(\epsilon,\rho)+\dim\ker
D_\epsilon(\infty)}{2}\right)\ .$$
It is now easy to see that $\tau(\epsilon,\rho)$ depends smoothly on the
parameters. In fact a change of the
kernel of  $D_\epsilon(\infty)$ gives rise to an integer jump of 
$\eta(\epsilon,\rho)$, and the combination $\eta(\epsilon,\rho)+\dim\ker
D_\epsilon(\infty)$ jumps by an even integer.

\begin{lem}\label{uin}
$\tau(\epsilon,\rho)$ is independent of $\rho$.
\end{lem}
\proof
The usual derivation for the variation of the $\eta$-invariant gives
$$\frac{\delta}{\delta\rho} \log \tau(\epsilon,\rho)=2\pi\imath \:\LIM_{t\to
0} \frac{1}{\pi^{1/2}}\Tr\:\left(\frac{\delta}{\delta\rho}
D_\epsilon(t)\right)\ee^{-D_\epsilon(t)^2} =0$$
since $\frac{\delta}{\delta\rho}
D_\epsilon(t)\equiv 0$ for small $t$.
\hB

From now on we assume that $\rho(s)\equiv 1$ for large $s$.
We form the family of super connections
$$A_s(\rho):=s \sigma D(\cE) + s \rho(s) \sigma Q + \nabla^{\Gamma(\cE)} +
\frac{\sigma}{4s} c(T)$$
and define the eta form
\begin{equation}\label{def1}\hat  \eta(\rho):=\cR\:
\frac{1}{\pi^{1/2}}\int_0^\infty \Tr_\sigma
\left(\frac{dA_s(\rho)}{ds} \ee^{-A^2_s(\rho)} \right) ds\ .\end{equation}
For small  $s$ the super connection $A_s(\rho)=A_s(\cE)$
 is the Bismut super connection and therefore the integral converges at $s=0$.
For large $s$ the operator  $D(\cE) + \rho(s)  Q=D(\cE)+Q$ is
invertible, and therfore the integrand vanishes exponentially for $s\to
\infty$.

The following proposition is an analog of the result of Bismut-Cheeger 
\cite{bismutcheeger89}.
\begin{prop} \label{limit}
We have $$\lim_{\epsilon\to 0}
\tau(\epsilon,\rho)=\exp\left( 2\pi\imath \int_{S^1\times S^1}\hat
\eta(\rho)\right)\ .$$ \end{prop}
We will prove this proposition in Subsection \ref{limitproof}.

If we replace the definition (\ref{def1}) of $\hat \eta(\rho)$ by
$$
\hat \eta(\rho):=\cR\: \LIM_{t\to 0}\frac{1}{\pi^{1/2}} \int_t^\infty
\Tr_\sigma
\left(\frac{dA_s(\rho)}{ds} \ee^{-A^2_s(\rho)} \right) ds\ ,$$
then it extends to all functions $\rho\in C^\infty[0,\infty)$
satisfying $\rho(s)\equiv 1$ for large $s$ and $\rho(s)\equiv const$ for small
$s$ (this is completely parallel to the definition  of
$\hat \eta_\alpha$ in \cite{lott01}, (2.19)). 

\begin{lem}\label{exact}
The difference $\hat \eta(\rho)-\hat\eta$ is exact.
\end{lem}
\proof
Let $\rho_u$, $u\in[0,1]$ be a smooth family of functions as above
interpolating between $\rho=\rho_0$ and the constant function $\rho_1\equiv 1$.
We now repeat the argument of the proof of \cite{lott01}, Proposition  3.
Fortunately, Lott has written the argument  in a way which does not use the
precise $s$-dependence of $A_s(\rho)$. In order to get \cite{lott01},
equation (2.11), one only needs that $\rho(s)$ is constant for small $s$.
\hB  
Combining Lemma \ref{exact} with Proposition \ref{limit}, and then further with
Lemma \ref{uin} and Equation (\ref{zwi})  we get 
\begin{eqnarray*} \exp\left(2\pi\imath \int_{S^1\times S^1} F \right)&=&
\exp\left(2\pi\imath \int_{S^1\times S^1}\hat  \eta(\rho) \right)\\
&=&\lim_{\epsilon\to 0}
\tau(\epsilon,\rho)\\&=&
\lim_{\epsilon\to
0}\tau(D_\epsilon)\\
&=&
T(\det(L_0\cE))(\Gamma)\ .
\end{eqnarray*}
In view of Equation (\ref{zw2}) this shows Proposition \ref{k1t}. \hB

\subsection{Canceling the spectral flow}

In this Subsection we combine Propositions \ref{bounding} and \ref{k1t}
in order to deduce Proposition \ref{stmain}.

Let $\cE_1$ be a geometric family over $S^1$ with odd dimensional fibres
such that $\ind(\cE_1)$ is the generator of $K^1(S^1)\cong \Z$. In fact, for
$\cE_1$ one can take the family given by $\pr:S^1\times S^1\rightarrow
S^1$ with standard metrics and horizontal subspaces such that $\cL$ is a line
bundle with connection over $S^1\times S^1$ with first Chern class satisfying
$c_1(\cL)[S^1\times S^1]=1$.

If we choose generators $\gamma_1,\gamma_2$ of $\pi_1(S^1\times S^1)\cong
\Z^2$, then we obtain an identification $K^1(S^1\times S^1)\cong \Z^2$. 
If $\cE$ is a geometric family over $S^1\times S^1$ with
odd-dimensional fibres, then $\ind(\cE)=(sf(\cE)(\gamma_1),sf(\cE)(\gamma_2))$,
where $sf(\cE)(\gamma)$ denotes the spectral flow along $\gamma$.

Let now $\cE$ be any geometric family with odd-dimensional fibres over some
base $B$. The spectral flow provides a homomorphism $\pi_1(B)\rightarrow \Z$.
This homomorphism corresponds to an element $sf(\cE)\in H^1(M,\Z)$.
Let $S:B\rightarrow S^1$ be the classifiying map of $-sf(\cE)$.
Then the spectral flow of $S^*\cE_1\cup \cE$ vanishes. In particular,
if $\Gamma\in L^2B$, then $\ind(\Gamma^*(S^*\cE_1 \cup_B\cE))=0$.

Therefore we can apply Proposition  \ref{k1t} in order to conclude
 $$T^2(\gerbe(S^*\cE_1 \cup_B\cE))(\Gamma)=T(\det(L_0(S^*\cE_1
\cup_B\cE)))(\Gamma)\ .$$
Using multiplicativity of gerbes and determinant bundles under relative
disjoint sum we obtain
\begin{equation}\label{zt1}\frac{T^2(\gerbe(\cE))(\Gamma)}{T(\det(L_0(\cE)))(\Gamma)}=\frac{T^2(\gerbe(S^*\cE_1))(\Gamma)}{T(\det(L_0(S^*\cE_1)))(\Gamma)} \ .\end{equation}
Now the map $S\circ \Gamma:S^1\times S^1\rightarrow S^1$ bounds. In fact, let
$u\in \pi_1(S^1\times S^1)$ be a primitive generator of $\ker  (S\circ
\Gamma)_*:\pi_1(S^1\times S^1)\rightarrow \pi_1(S^1)$. Then can find a
diffeomorphism  $F:S^1\times S^1\rightarrow S^1\times S^1$ such that
$F(S^1\times\{*\})$  represents $u$. The map $S\circ \Gamma \circ F$
now extends to $D^2\times S^1$. It follows that we can apply Proposition 
\ref{bounding} in order to conclude that the right hand side of (\ref{zt1}) is
equal to one. This proves Proposition  \ref{stmain} in general (up to the
verification of Proposition \ref{limit}). \hB

\subsection{Adiabatic limit of $\eta$-invariants}\label{limitproof}

In this subsection we prove Proposition  \ref{limit}.
In fact we can repeat the proof of Bismut-Cheeger \cite{bismutcheeger89}, Thm.
4.35. We will just explain the changes in various places. 
The main point is that we
must replace $u^{1/2}D_\epsilon$ by
$u^{1/2}(D_\epsilon(\cE)+\rho(u^{1/2})Q)=:D_\epsilon(u^{1/2})$. The right hand
side of formula \cite{bismutcheeger89}, (4.39),  gets correspondingly replaced
by  $$\frac{1}{2}\left[\dim\ker
D_\epsilon(\infty)+\frac{2}{\pi^{1/2}}\int_0^\infty \Tr\:
\left(\frac{\partial}{\partial u}
D_\epsilon(u^{1/2})\right)\ee^{-D^2_\epsilon(u^{1/2})}du\right]\ .$$ The proof
of \cite{bismutcheeger89}, Proposition  4.41, can be applied to $D_\epsilon(u^{1/2})$
for $u\ge 4$ and yields a lower bound $\spec(D^2_\epsilon(u^{1/2}))\ge
\lambda_0>0$ which is uniform for $u\ge 4$ and $\epsilon\in(0,\epsilon_0]$,
$\epsilon_0$ sufficiently small.  

In step ii of the proof of  \cite{bismutcheeger89}, Thm.
4.35 we must replace, of course,
$u D_\epsilon^2 - z u^{1/2} D_\epsilon$ by
$D_\epsilon(u^{1/2})^2-z 2u \frac{\partial}{\partial u} D_\epsilon(u^{1/2})$.
From  \cite{bismutcheeger89}, (4.45), we get (note that $D^Z$ in
\cite{bismutcheeger89} is $D(\cE)$ here)
 \begin{eqnarray*}
D_\epsilon(u^{1/2})- u^{1/2} D_\epsilon&=&u^{1/2}\rho(u^{1/2})Q\\
 D_\epsilon(u^{1/2})^2- uD^2_\epsilon&=&\epsilon^{1/2}u^{1/2}\rho(u^{1/2})
[E_\epsilon,Q] + u\rho(u^{1/2})  \{D(\cE), Q\} + u\rho(u^{1/2})^2 Q^2\\
2u\frac{\partial}{\partial u} D_\epsilon(u^{1/2})- u^{1/2}
D_\epsilon&=& u\rho^\prime(u^{1/2}) Q
\end{eqnarray*}
Therefore, in the right hand side of \cite{bismutcheeger89},
(4.51), we have the following additional terms
\begin{eqnarray*}
 &&\epsilon^{1/2}u^{1/2}\rho(u^{1/2})
[E_\epsilon,Q] + u\rho(u^{1/2})  \{D(\cE), Q\} + u\rho(u^{1/2})^2 Q^2\\
&& -z u\rho^\prime(u^{1/2}) Q\ .
\end{eqnarray*}
We now proceed as in 
\cite{bismutcheeger89} and  apply Getzler's rescaling along the base.
This is possible since the perturbation $Q$ is local with respect to the base.
Formula \cite{bismutcheeger89}, (4.70), gets replaced by 
$$\cH+A^2_{u^{1/2}}(\rho) - z 2u \frac{\partial}{\partial u} A_{u^{1/2}}(\rho)\
$$
(after setting $\sigma=1$).
The Hermite operator $\cH$ simplifies to $-\sum_\alpha(\partial_\alpha)^2$
since the torus $S^1\times S^1$ is flat.
Equation \cite{bismutcheeger89} (4.74) simplifies to 
\begin{equation}\label{yt1}\frac{1}{\pi^{1/2}} \Tr\:
\left(\frac{\partial}{\partial u} D_\epsilon(u)\right)
\ee^{-D_\epsilon^2(u)}\stackrel{\epsilon\to 0}{\longrightarrow}
\frac{1}{2\pi\imath}\frac{1}{\pi^{1/2}}\int_{S^1\times S^1}\Tr_\sigma
\frac{\partial}{\partial u} A_{u^{1/2}}(\rho) \ee^{-A^2_{u^{1/2}}(\rho)}  \
.\end{equation} Since $\rho(s)$ vanishes for small $s$, the discussion in
\cite{bismutcheeger89} of uniformity of the convergence (\ref{yt1}) on
intervals $(0,T]$, $T>0$ applies to the present case. Again using the fact,
that $Q$ is local with respect to the base we repeat the derivation of the
remainder term $O(\epsilon^{1/2}(1+T^N))$ in (\ref{yt1}). Thus we have
verified the analog of \cite{bismutcheeger89} (4.40). The final step iii. of
the proof of \cite{bismutcheeger89}, Thm. 4.35, can be taken without change.
\hB

\bibliographystyle{plain}
%\bibliography{/home/uli/working/libank/literatu,/home/uli/working/libank/lit1} 
%\bibliography{/user/bunke/working/libank/literatu,/user/bunke/working/libank/lit1}

\end{document}